\documentclass{amsart}%
\usepackage{amssymb}
\usepackage{amsfonts}
\usepackage{amsmath}
\usepackage{graphicx}%
\newtheorem{theorem}{Theorem}
\theoremstyle{plain}

\begin{document}
\title[Remark on a geometric inequality]{Remark on an inequality for closed hypersurfaces in complete manifolds with
nonnegative Ricci curvature}
\author{Xiaodong Wang}
\address{Department of Mathematics, Michigan State University, East Lansing, MI 48824}
\email{xwang@math.msu.edu}
\maketitle

\begin{abstract}
We give a simple proof of a recent result due to Agostiniani, Fogagnolo and
Mazzieri \cite{AFM}.

\end{abstract}

The following result was proved by Agostiniani, Fogagnolo and Mazzieri
\cite{AFM}.

\begin{theorem}
Let $\left(  M^{n},g\right)  $ ($n\geq3$) be a complete Riemannian manifold
with nonnegative Ricci curvature and $\Omega\subset M$ a bounded open set with
smooth boundary. Then
\begin{equation}
\int_{\partial\Omega}\left\vert \frac{H}{n-1}\right\vert ^{n-1}d\sigma
\geq\mathrm{AVR}\left(  g\right)  \left\vert \mathbb{S}^{n-1}\right\vert ,
\label{Will}%
\end{equation}
where $H$ is the mean curvature of $\partial\Omega$ and $\mathrm{AVR}\left(
g\right)  $ is the asymptotic volume ratio of $M$. Moreover, if $\mathrm{AVR}%
\left(  g\right)  >0$, equality holds iff $M\backslash\Omega$ is isometric to
$\left(  [r_{0},\infty)\times\partial\Omega,dr^{2}+\left(  \frac{r}{r_{0}%
}\right)  ^{2}g_{\partial\Omega}\right)  $ with
\[
r_{0}=\left(  \frac{\left\vert \partial\Omega\right\vert }{\mathrm{AVR}\left(
g\right)  \left\vert \mathbb{S}^{n-1}\right\vert }\right)  ^{\frac{1}{n-1}}%
\]
In particular, $\partial\Omega$ is a connected totally umbilic submanifold
with constant mean curvature.
\end{theorem}

\bigskip The proof in \cite{AFM} is highly nontrivial. It is based on the
study of the solution of the following problem
\[
\left\{
\begin{array}
[c]{cc}%
\Delta u=0, & \text{on }M\backslash\Omega\\
u=1 & \text{on \ \ }\partial\Omega\\
u\left(  x\right)  \rightarrow0 & \text{as }x\rightarrow\infty,
\end{array}
\right.
\]
which exists when $\mathrm{AVR}\left(  g\right)  >0$. The key step consists of
showing that, with $\beta\geq\left(  n-2\right)  /\left(  n-1\right)  $%
\[
U_{\beta}\left(  t\right)  =t^{-\beta\left(  \frac{n-1}{n-2}\right)  }%
\int_{u=t}\left\vert \nabla u\right\vert ^{\beta+1}d\sigma
\]
is monotone in $t\in(0,1]$. The geometric inequality (\ref{Will}) then follows
by analyzing the asymptotic behavior of $U_{\beta}\left(  t\right)  $ as
$t\rightarrow0$. It is a beautiful argument.

In this short note, we show that this theorem can be proved by standard
comparison methods in Riemannian geometry.

To prove the inequality (\ref{Will}), we assume, without loss of generality,
that $\Omega$ has no hole, i.e. $M\backslash\Omega$ has no bounded component.
In the following we write $\Sigma=\partial\Omega$ and let $\nu$ be the outer
unit normal along $\Sigma$. For each $p\in\Sigma$ let $\gamma_{p}\left(
t\right)  =\exp_{p}t\nu\left(  p\right)  $ be the normal geodesic with initial
velocity $\nu\left(  p\right)  $. We define
\[
\tau\left(  p\right)  =\sup\left\{  L>0:\gamma_{p}\text{ is minimizing on
}\left[  0,L\right]  \right\}  \in(0,\infty].
\]
It is well known that $\tau$ is a continuous function on $\Sigma$ and the
focus locus
\[
C\left(  \Sigma\right)  =\left\{  \exp_{p}\tau\left(  p\right)  \nu\left(
p\right)  :\tau\left(  p\right)  <\infty\right\}
\]
is a closed set of measure zero in $M$. Moreover the map $\Phi\left(
r,p\right)  =\exp_{p}r\nu\left(  p\right)  $ is a diffeomorphism from
\[
E=\left\{  \left(  r,p\right)  \in\Sigma\times\lbrack0,\infty):r<\tau\left(
p\right)  \right\}
\]
onto $\left(  M\backslash\Omega\right)  \backslash C\left(  \Sigma\right)  $.
And on $E$ the pull back of the volume form takes the form $d\mu
=\mathcal{A}\left(  r,p\right)  drd\sigma\left(  p\right)  $. We will also
understand $r$ as the distance function to $\Sigma$ and it is smooth on
$M\backslash\Omega$ away from $C\left(  \Sigma\right)  $. By the Bochner
formula and nonnegative Ricci curvature condition%
\begin{align*}
0  &  =\frac{1}{2}\Delta\left\vert \nabla r\right\vert ^{2}=\left\vert
D^{2}r\right\vert ^{2}+\left\langle \nabla r,\nabla\Delta r\right\rangle
+Ric\left(  \nabla r,\nabla r\right) \\
&  \geq\frac{\left(  \Delta r\right)  ^{2}}{n-1}+\frac{\partial}{\partial
r}\Delta r.
\end{align*}
In view of the initial condition $\Delta r|_{r=0}=H$, it is standard to deduce
from the above inequality $\tau\leq\frac{n-1}{H^{-}}$ and
\[
\frac{\mathcal{A}^{\prime}}{\mathcal{A}}=\Delta r\leq\frac{\left(  n-1\right)
H}{n-1+Hr}%
\]
This shows that the function%
\[
\theta\left(  r,p\right)  =\frac{\mathcal{A}\left(  r,p\right)  }{\left(
1+\frac{H\left(  p\right)  }{n-1}r\right)  ^{n-1}}%
\]
is non-increasing in $r$ on $[0,\tau\left(  p\right)  )$. As $\theta\left(
0,p\right)  =1$, we obtain
\[
\mathcal{A}\left(  r,p\right)  \leq\left(  1+\frac{H\left(  p\right)  }%
{n-1}r\right)  ^{n-1}.
\]
The above analysis is standard in Riemannian geometry. We also remark that
this argument involving the Bochner formula can be replaced by an argument
involving the index form along each individual geodesic $\gamma_{p}$. For more
details, cf. \cite{P, S} or other books on Riemannian geometry.

Therefore for any $R>0$%
\begin{align*}
Vol\left\{  x\in M:d\left(  x,\Omega\right)  <R\right\}  =  &  \left\vert
\Omega\right\vert +\int_{\Sigma}\int_{0}^{\min\left(  R,\tau\left(  p\right)
\right)  }\mathcal{A}\left(  r,p\right)  drd\sigma\left(  p\right) \\
\leq &  \left\vert \Omega\right\vert +\int_{\Sigma}\int_{0}^{\min\left(
R,\tau\left(  p\right)  \right)  }\left(  1+\frac{H\left(  p\right)  }%
{n-1}r\right)  ^{n-1}drd\sigma\left(  p\right) \\
\leq &  \left\vert \Omega\right\vert +\int_{\Sigma}\int_{0}^{\min\left(
R,\tau\left(  p\right)  \right)  }\left(  1+\frac{H^{+}\left(  p\right)
}{n-1}r\right)  ^{n-1}drd\sigma\left(  p\right) \\
\leq &  \left\vert \Omega\right\vert +\int_{\Sigma}\int_{0}^{R}\left(
1+\frac{H^{+}\left(  p\right)  }{n-1}r\right)  ^{n-1}drd\sigma\left(  p\right)
\\
=  &  \left\vert \Omega\right\vert +\frac{R^{n}}{n}\int_{\Sigma}\left(
\frac{H^{+}\left(  p\right)  }{n-1}\right)  ^{n-1}d\sigma\left(  p\right)
+O\left(  R^{n-1}\right)  .
\end{align*}
Dividing both sides by $\left\vert \mathbb{B}^{n}\right\vert R^{n}=\left\vert
\mathbb{S}^{n-1}\right\vert R^{n}/n$ and letting $R\rightarrow\infty$ yields
\[
\mathrm{AVR}\left(  g\right)  \leq\frac{1}{\left\vert \mathbb{S}%
^{n-1}\right\vert }\int_{\Sigma}\left(  \frac{H^{+}}{n-1}\right)
^{n-1}d\sigma,
\]
which implies (\ref{Will}).

We now analyze the equality case. Suppose
\begin{equation}
\mathrm{AVR}\left(  g\right)  =\frac{1}{\left\vert \mathbb{S}^{n-1}\right\vert
}\int_{\Sigma}\left(  \frac{H^{+}}{n-1}\right)  ^{n-1}d\sigma>0.\label{SWill}%
\end{equation}
It is clear from the proof that $\tau\equiv\infty$ on the open set $\Sigma
^{+}=\left\{  p\in\Sigma:H\left(  p\right)  >0\right\}  $. For any $R^{\prime
}<R$ we have%
\begin{align*}
&  Vol\left\{  x\in M:d\left(  x,\Omega\right)  <R\right\}  \\
= &  \left\vert \Omega\right\vert +\int_{\Sigma^{+}}\int_{0}^{R}%
\mathcal{A}\left(  r,p\right)  drd\sigma\left(  p\right)  +\int_{\Sigma
\backslash\Sigma^{+}}\int_{0}^{\min\left(  R,\tau\left(  p\right)  \right)
}\mathcal{A}\left(  r,p\right)  drd\sigma\left(  p\right)  \\
\leq &  \left\vert \Omega\right\vert +\int_{\Sigma^{+}}\int_{0}^{R}%
\theta\left(  r,p\right)  \left(  1+\frac{H\left(  p\right)  }{n-1}r\right)
^{n-1}drd\sigma\left(  p\right)  +\int_{\Sigma\backslash\Sigma^{+}}\int
_{0}^{R}drd\sigma\left(  p\right)  \\
\leq &  \left\vert \Omega\right\vert +\int_{\Sigma^{+}}\int_{R^{\prime}}%
^{R}\theta\left(  r,p\right)  \left(  1+\frac{H\left(  p\right)  }%
{n-1}r\right)  ^{n-1}drd\sigma\left(  p\right)  \\
&  +\int_{\Sigma^{+}}\int_{0}^{R^{\prime}}\theta\left(  r,p\right)  \left(
1+\frac{H\left(  p\right)  }{n-1}r\right)  ^{n-1}drd\sigma\left(  p\right)
+O\left(  R\right)  \\
\leq &  \left\vert \Omega\right\vert +\int_{\Sigma^{+}}\theta\left(
R^{\prime},p\right)  \int_{R^{\prime}}^{R}\left(  1+\frac{H\left(  p\right)
}{n-1}r\right)  ^{n-1}drd\sigma\left(  p\right)  \\
&  +\int_{\Sigma^{+}}\int_{0}^{R^{\prime}}\theta\left(  r,p\right)  \left(
1+\frac{H\left(  p\right)  }{n-1}r\right)  ^{n-1}drd\sigma\left(  p\right)
+O\left(  R\right)  .
\end{align*}
Dividing both sides by $\left\vert \mathbb{B}^{n}\right\vert R^{n}=\left\vert
\mathbb{S}^{n-1}\right\vert R^{n}/n$ and letting $R\rightarrow\infty$ yields%
\[
\mathrm{AVR}\left(  g\right)  \leq\frac{1}{\left\vert \mathbb{S}%
^{n-1}\right\vert }\int_{\Sigma^{+}}\left(  \frac{H\left(  p\right)  }%
{n-1}\right)  ^{n-1}\theta\left(  R^{\prime},p\right)  d\sigma\left(
p\right)  .
\]
Letting $R^{\prime}\rightarrow\infty$ yields%
\[
\mathrm{AVR}\left(  g\right)  \leq\frac{1}{\left\vert \mathbb{S}%
^{n-1}\right\vert }\int_{\Sigma^{+}}\left(  \frac{H}{n-1}\right)  ^{n-1}%
\theta_{\infty}d\sigma,
\]
where $\theta_{\infty}\left(  p\right)  =\lim_{r\rightarrow\infty}%
\theta\left(  r,p\right)  \leq1$. As we have equality (\ref{SWill}) we must
have $\theta_{\infty}\left(  p\right)  =1$ for a.e. $p\in\Sigma^{+}$. It
follows that
\[
\mathcal{A}\left(  r,p\right)  =\left(  1+\frac{H\left(  p\right)  }%
{n-1}r\right)  ^{n-1}\text{ on }[0,\infty)
\]
for a.e. $p\in\Sigma^{+}$. By continuity the above identity holds for all
$p\in\Sigma^{+}$.

Inspecting the comparison argument, we must have on $\Phi\left(
\lbrack0,\infty)\times\Sigma^{+}\right)  $
\begin{align*}
D^{2}r  &  =\frac{\Delta r}{n-1}g=\frac{H}{n-1+Hr}g,\\
Ric\left(  \nabla r,\nabla r\right)   &  =0.
\end{align*}
As $Ric\geq0$, it follows that $Ric\left(  \nabla r,\cdot\right)  =0$. From
the 1st equation above $\Sigma^{+}$ is an umbilic hypersurface, i.e. the 2nd
fundamental form $\Pi=\frac{H}{n-1}g_{\Sigma^{+}}$. Working with an
orthonormal frame $\left\{  e_{0}=\nu,e_{1},\cdots,e_{n-1}\right\}  $ along
$\Sigma^{+}$ we have by the Codazzi equation, with $1\leq i,j,k\leq n-1$%
\[
R\left(  e_{k},e_{j},e_{i},\nu\right)  =\Pi_{ij,k}-\Pi_{ik,j}=\frac{1}%
{n-1}\left(  H_{k}\delta_{ij}-H_{j}\delta_{ik}\right)  .
\]
Taking trace over $i$ and $k$ yields
\[
-\frac{n-2}{n-1}H_{j}=Ric\left(  e_{j},\nu\right)  =0.
\]
As a result $H$ is locally constant on $\Sigma^{+}$. Therefore $\Sigma^{+}$
must be the union of several components of $\Sigma$. We know that $\Phi$ is a
diffeomorphism form $[0,\infty)\times\Sigma^{+}$ onto its image and the
pullback metric $\Phi^{\ast}g$ takes the following form%
\[
dr^{2}+h_{r},
\]
where $h_{r}$ is a $r$-dependent family of metrics on $\Sigma^{+}$ and
$h_{0}=g_{\Sigma^{+}}$. We have
\[
D^{2}r=\frac{H}{n-1+Hr}g.
\]
In terms of local coordinates $\left\{  x_{1},\cdots,x_{n-1}\right\}  $ on
$\Sigma^{+}$ the above equation implies%
\[
\frac{1}{2}\frac{\partial}{\partial r}h_{ij}=\frac{H}{n-1+Hr}h_{ij}.
\]
Therefore $h_{r}=\left(  1+\frac{H}{n-1}r\right)  ^{2}g_{\Sigma^{+}}$. This
proves that $\Phi\left(  \lbrack0,\infty)\times\Sigma^{+}\right)  $ is
isometric to $\left(  [r_{0},\infty)\times\Sigma^{+},dr^{2}+\left(  \frac
{r}{r_{0}}\right)  ^{2}g_{\Sigma^{+}}\right)  $, where $r_{0}=\frac{n-1}{H}$.

Since $M$ has nonnegative Ricci curvature and Euclidean volume growth, it has
only one end by the Cheeger-Gromoll theorem. Therefore $\Sigma^{+}$ is
connected and if $\Sigma$ has other components besides $\Sigma^{+}$, they all
bound bounded components of $M\backslash\Omega$.

If we have the stronger identity%
\[
\mathrm{AVR}\left(  g\right)  =\frac{1}{\left\vert \mathbb{S}^{n-1}\right\vert
}\int_{\Sigma}\left\vert \frac{H}{n-1}\right\vert ^{n-1}d\sigma>0,
\]
inspecting the proof of the inequality (\ref{Will}) shows that we must have
$H\geq0$ on $\Sigma$. Then $\overline{\Omega}$ is compact Riemannian manifold
with mean convex boundary. It is a classic fact that its boundary must be
connected, see \cite{I, K} or \cite{HW} for an analytic argument. Therefore
$\Sigma=\Sigma^{+}$ is connected and $M\backslash\Omega$ is isometric to
$\left(  [r_{0},\infty)\times\Sigma,dr^{2}+\left(  \frac{r}{r_{0}}\right)
^{2}g_{\Sigma}\right)  $.

\bigskip

\textbf{Acknowledgement}.\textbf{ }I would like to thank Fengbo Hang for
helpful discussions.


\begin{thebibliography}{999}                                                                                              %


\bibitem[AFM]{AFM}V. Agostiniani; M. Fogagnolo; L. Mazzieri. Sharp geometric
inequalities for closed hypersurfaces in manifolds with nonnegative Ricci
curvature. Invent. Math. 222 (2020), no. 3, 1033-1101.

\bibitem[HW]{HW}F. Hang; X. Wang. Vanishing sectional curvature on the
boundary and a conjecture of Schroeder and Strake, Pacific J. Math. 232
(2007), no. 2, 283-287.

\bibitem[I]{I}R. Ichida. Riemannian manifolds with compact boundary. Yokohama
Math. J. 29 (1981), no. 2, 169-177.

\bibitem[K]{K}A. Kasue. Ricci curvature, geodesics and some geometric
properties of Riemannian manifolds with boundary. J. Math. Soc. Japan 35
(1983), no. 1, 117-131.

\bibitem[P]{P}P. Petersen. Riemannian Geometry. Third edition. Graduate Texts
in Mathematics, 171. Springer, Cham, 2016.

\bibitem[S]{S}T. Sakai. Riemannian geometry. Translations of Mathematical
Monographs, 149. American Mathematical Society, Providence, RI, 1996.
\end{thebibliography}
\end{document}